\documentclass[journal ]{new-aiaa}
\usepackage[utf8]{inputenc}
\usepackage{textcomp}
\usepackage{ulem} 
\usepackage{graphicx}
\usepackage{amsmath}
\usepackage[version=4]{mhchem}
\usepackage{siunitx}
\usepackage{longtable,tabularx}
\usepackage{color}
\usepackage{mathtools}
\usepackage{float}
\usepackage{caption,subcaption}
\usepackage{tikz}
\usepackage[many]{tcolorbox}
\usepackage{color}
\usepackage{epstopdf}
\usepackage[numbers,sort&compress]{natbib}
\usepackage{algorithm}
\usepackage{algorithmic}
\usepackage{setspace}

\def\1{\boldsymbol{1}}

\title{Fuel-Optimal Trajectory Planning for Lunar Vertical Landing}

\author{Kun Wang\footnote{PhD student, School of Aeronautics and Astronautics.}, Zheng Chen\footnote{Researcher, School of Aeronautics and Astronautics. email: \underline{z-chen@zju.edu.cn}.}, and Jun Li\footnote{Professor, School of Aeronautics and Astronautics.}}
\affil{ Zhejiang University, Hangzhou 310027, Zhejiang, China}
\begin{document}
\maketitle

\begin{abstract}
In this paper, we consider a trajectory planning problem arising from a lunar vertical landing with minimum fuel consumption. The vertical landing requirement is written as a final steering angle constraint, and a nonnegative regularization term is proposed to modify the cost functional. In this way, the final steering angle constraint will be inherently satisfied according to Pontryagin's Minimum Principle. As a result, the modified optimal steering angle has to be determined by solving a transcendental equation. To this end, a transforming procedure is employed, which allows for finding the desired optimal steering angle by a simple bisection method. Consequently, the modified optimal control problem can be solved by the indirect shooting method. Finally, some numerical examples are presented to demonstrate and verify the developments of the paper. 
\end{abstract}

\section{ Introduction}
\label{intro}
As the nearest celestial body to the Earth, exploration of the Moon has been receiving a great deal of attention recently. Soft landing is the key to deploying the lunar lander for scientific missions. The typical soft landing consists of two phases, i.e., de-orbit burn phase and powered descent phase \cite{leeghim2016optimal}. During the de-orbit burn phase, the altitude of the lunar lander is reduced from the parking orbit to the perilune altitude. On the other hand, the powered descent phase is expected to use the retro-propulsion provided by the engine to meet the required final condition \cite{lu2018propellant}. In addition, minimizing the fuel consumption is also of great importance. Thus, significant efforts have been devoted to designing the fuel-optimal landing trajectory, which is equivalent to solving an Optimal Control Problem (OCP).

An OCP can be numerically solved via the direct or the indirect shooting method, as summarized in Ref. \cite{betts1998survey}. The direct methods, especially convex optimization based methods, have been seeing wide use in onboard implementation for soft landing problems; see, e.g., \cite{accikmecse2013lossless,sagliano2023six}. On the other hand, the indirect shooting method employs Pontryagin's Minimum Principle (PMP) to transform the OCP into a Two-Point Boundary-Value Problem (TPBVP). Albeit not ideal for onboard implementation, this method has been widely used to reveal the structure of the optimal control law; see, e.g., \cite{lu2018propellant,ito2020throttled,lu2023propellant}. Moreover, learning-based methods have been employed to generate the optimal or suboptimal solution in real time, where neural networks are trained to approximate the optimal control; see, e.g., \cite{gaudet2014adaptive,gaudet2020deep,you2021learning,cheng2020real,wang2022nonlinear,wang2023real}. 

Most of these methods mentioned in the preceding paragraph, whether the obtained solution is real-time or not, have solely focused on OCPs without considering the final control constraint. 
From a practical point of view, the final control constraint can be imposed to guarantee smooth system operation \cite{roozegar2018optimal}, reduce final errors \cite{lee2003optimal}, or improve the missile's guidance performance \cite{ryoo2006time,lee2013polynomial,qiu2014research,WangZhi2021}, to name a few. For the lunar landing problem, the thrust engine is {\it de facto} usually fixed with the lunar lander's body. Therefore, the steering angle essentially represents the attitude of the lunar lander \cite{lu2023propellant}. A vertical landing, characterized by a final steering angle constraint, can prevent the lunar lander from rollover, or be further used for detecting obstructions and implementing terminal avoidance manoeuvres \cite{mcinnes1995path}. For this reason, there were some works to generate the fuel-optimal vertical landing trajectory; see, e.g., \cite{mcinnes1995path,zhou2010optimal,sachan2015fuel,leeghim2016optimal,lu2018propellant,wang2024neural}. Specifically, in Ref. \cite{leeghim2016optimal}, the lunar landing trajectory was deconstructed into two phases, and the final steering angle constraint was augmented into the cost functional along the second phase. Then, the original TPBVP became a Multi-PBVP. In Ref. \cite{lu2018propellant}, a better estimation of the time-to-go from Ref. \cite{citron1964terminal} was implemented into the Apollo guidance proposed in Refs. \cite{cherry1964general,klumpp1974apollo}, and the vertical landing was realized by using extra fuel.  In Ref. \cite{zhou2010optimal}, a control parameterization and time scaling transform technique were used to realize a soft landing with the terminal attitude of the lunar lander to be within a small deviation from vertical landing. In Ref. \cite{sachan2015fuel}, to deal with the final steering angle constraint,  a time varying matrix was introduced to the cost functional, and the resulting problem was solved via  model predictive static programming. However, the non-analytical time varying matrix may not be generalizable to broader cases. In Ref. \cite{wang2024neural}, a nonnegative regularization term was proposed to augment the cost functional, ensuring that the final steering angle constraint is met. Nevertheless, in that work, the thrust magnitude of the lunar lander is not adjustable.

In this paper, in order to cater to the final steering angle constraint, which can be seen as a final control constraint, we modify the cost functional by augmenting a simple nonnegative regularization term, as done in Ref. \cite{wang2024neural}. In this way, the final steering angle constraint will be satisfied at the final time automatically because Pontryagin's Minimum Principle (PMP) minimizes the augmented cost functional. The proposed method can be readily embedded into the indirect shooting method once the modified optimal steering angle is obtained by solving a transcendental equation, whose desired zero can be found by a simple bisection method. Finally, numerical simulations are presented, showing that the proposed method is able to find the fuel-optimal landing trajectory with only a negligible penalty on the fuel consumption, and the final steering angle constraint ensuring a vertical landing can be strictly satisfied.

The remainder of the paper is organized as follows. Section \ref{sec2} presents the OCP regarding the lunar vertical landing. PMP is used to derive the necessary conditions for optimality without and with the final steering angle constraint in Sections \ref{sec3} and \ref{sec4}, respectively. Numerical simulations are presented in Section \ref{sec5}, followed by Section \ref{sec6} for conclusions. 

\section{Problem Formulation}\label{sec2}
\begin{figure}[!htp]
\begin{center}
\includegraphics[scale=0.3]{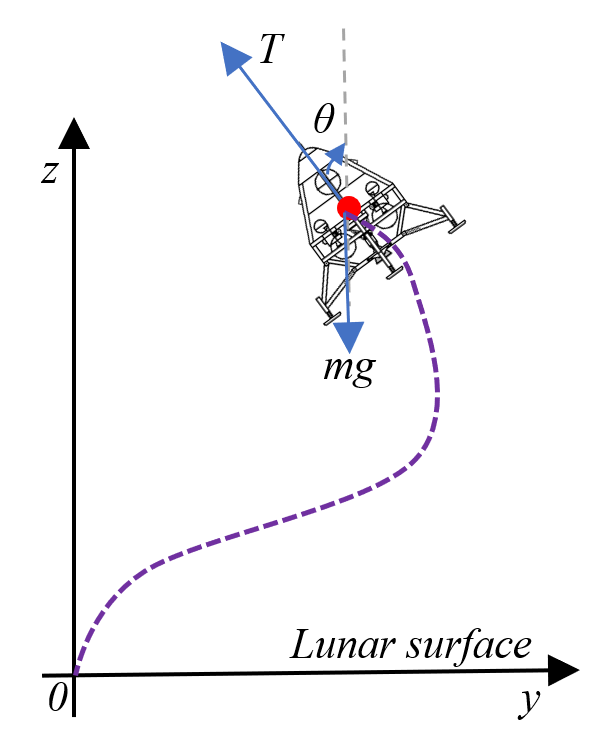}
\caption{Coordinate system for the lunar vertical landing.}\label{Fig:frame}
\end{center}
\end{figure}

Consider the 2-dimensional lunar soft landing problem in Ref. \cite{sanchez2018real}, as shown in Fig.~\ref{Fig:frame}, where a flat non-rotating Moon model with constant gravity is assumed. The origin of the coordinate system $Oyz$ is located at the landing site, with $y \in \mathbb{R}$ representing the lunar lander's ground range and $z \geq 0$ denoting its altitude.
The lunar lander is propelled by an engine with adjustable thrust magnitude and  steering angle. The thrust magnitude is denoted by $T \in [0,T_{m}]$, where the constant $T_{m}$ means the maximum thrust magnitude. Let $u\in[0,1]$ denote the engine thrust ratio, then we have $T=uT_{m}$. The steering angle, denoted by $\theta \in [-\pi, \pi]$, is defined to be the angle from the thrust vector to the local vertical line of the lunar lander. In addition, the gravity of the lunar lander is $mg$, in which $m$ is its mass, and $g$ is the Moon’s gravitational acceleration.
Denote by $\boldsymbol r = [y;z]$ the position vector, $\boldsymbol v = [v_y;v_z]$ the velocity vector of the lunar lander, respectively; $v_y$ is the horizontal speed and $v_z$ is the vertical speed.
Let $\boldsymbol x = [\boldsymbol r, \boldsymbol v, m]^T$ be the 
state vector.
Then, the dynamics of the lunar lander is represented by
\begin{align}
\dot{\boldsymbol x}(t) = \boldsymbol f(\boldsymbol x, u, \theta,t) \Rightarrow 
\begin{cases}
\dot{\boldsymbol r}(t) =  \boldsymbol v(t),\\
\dot{\boldsymbol v}(t) =  \boldsymbol g + \frac{u(t)T_{m}}{m(t)}\hat{\boldsymbol i}_\theta(t),\\
\dot{m}(t) = -\frac{u(t)T_{m}}{I_{sp}g_0},
\end{cases}
\label{EQ:dynamics}
\end{align}
where $t \geq 0$ is the time; $\boldsymbol f:\mathbb{R}^5 \times \mathbb{R} \times \mathbb{R} \times \mathbb{R}_0^+ \rightarrow \mathbb{R}^5$ is the smooth vector field. The acceleration vector $\boldsymbol g = [0;-g]$. $\hat{\boldsymbol i}_\theta = [\sin \theta;\cos \theta]$ is the unit thrust vector. The constant $I_{sp}$ represents the engine's specific impulse and $g_0$ is the Earth’s gravitational acceleration at sea level. 
The initial condition of the lunar lander at $t_0 = 0$ is given by
\begin{align}
\boldsymbol{r}(0) = \boldsymbol{r}_0, \boldsymbol{v}(0) = \boldsymbol{v}_0, m(0) = m_0.
\label{EQ:initialcondition}
\end{align}
The final condition at touchdown is 
\begin{align}
\boldsymbol{r}(t_f) = \boldsymbol{0}, \boldsymbol{v}(t_f) = \boldsymbol{0},
\label{EQ:finalcondition}
\end{align}
in which $t_f$ is the free final time. 
Because the steering angle defines the body attitude of the lunar lander \cite{lu2023propellant}, to ensure a vertical landing, the final steering angle should satisfy the following constraint:
\begin{align}
\theta(t_f) = 0.
\label{EQ:steeringcondition}
\end{align}
The cost functional to be minimized is
\begin{align}
J = \int_0^{t_f} u(t) ~\mathrm{d}t
\label{EQ:cost}
\end{align}

In the following two sections, we shall present the corresponding TPBVPs without and with the final steering angle constraint.
\section{Without the Final Steering Angle Constraint}\label{sec3}
Denote by $\boldsymbol{p}_r = [p_y;p_z]$, $\boldsymbol{p}_v=[p_{v_y};p_{v_z}]$, and $p_{m}$ the co-state variables of $\boldsymbol{r}$, $\boldsymbol{v}$, and $m$, respectively. Then, the Hamiltonian is expressed as
\begin{align}
\mathscr{H} = \boldsymbol{p}^T_r \cdot \boldsymbol{v}  + \boldsymbol{p}^T_v \cdot (\boldsymbol g + \frac{uT_{m}}{m}\hat{\boldsymbol i}_\theta) + p_{m} (-\frac{uT_{m}}{I_{sp}g_0}) + u
\label{EQ:Ham}
\end{align}
According to PMP \cite{Pontryagin}, it holds that
\begin{align}
\begin{cases}
\dot{\boldsymbol{p}}_r(t) = -\frac{\partial \mathscr{H}}{\partial \boldsymbol{r}(t)} = \boldsymbol{0},\\
\dot{\boldsymbol{p}}_v(t) = -\frac{\partial \mathscr{H}}{\partial \boldsymbol{v}(t)} = -\boldsymbol{p}_r(t),\\
\dot{p}_m(t) = -\frac{\partial \mathscr{H}}{\partial m(t)} =  \frac{u(t)T_m}{m^2(t)}\boldsymbol{p}^T_v(t) \cdot \hat{\boldsymbol i}_\theta(t).
\end{cases}
\label{EQ:dot_p}
\end{align}
The optimal steering angle follows that
\begin{align}
\hat{\boldsymbol i}^*_\theta(t) = -\frac{\boldsymbol{p}_v(t)}{\Vert \boldsymbol{p}_v(t)\Vert_2}.\label{EQ:dH/dtheta}
\end{align}
It is shown in Ref. \cite{lu2018propellant} that the optimal thrust magnitude is bang-bang without singular arc, thus the optimal thrust ratio is determined by
\begin{align}
u^*(t) = 
\begin{cases}
1, S(t) < 0, \\
0, S(t) > 0,
\end{cases}
\label{EQ:optimalration}
\end{align}
in which $S(t)$ is the switching function satisfying
\begin{align}
S(t) = 1 - \frac{T_m~p_m(t)}{I_{sp}g_0} - \frac{T_m}{m(t)}\Vert \boldsymbol{p}_v(t)\Vert_2.
\label{EQ:SF}
\end{align}
To avoid numerical difficulties caused by discontinuous bang-bang control, the smoothing technique in Ref. \cite{wang2023new} is adopted to approximate Eq.~(\ref{EQ:optimalration}), i.e.,
\begin{align}
u^*(t) \approx u^*(t,\delta) = \frac{1}{2}(1 - \frac{S(t)}{\sqrt{\delta + |S(t)|^2}}),
\label{EQ:smoothcontrol}
\end{align}
in which $\delta$ is the smoothing constant.
Since the final mass is free, the transversality condition implies
\begin{align}
p_{m}(t_f) = 0.
\label{EQ:Transversality1}
\end{align}
In addition, since the Hamiltonian in Eq.~(\ref{EQ:Ham}) does not contain time explicitly and the final time is free, the stationary condition states that
\begin{align}
\mathscr H(t_f) = 0.
\label{EQ:stat1}
\end{align}

Eqs.~(\ref{EQ:dynamics}), (\ref{EQ:dot_p}), (\ref{EQ:dH/dtheta}) and (\ref{EQ:smoothcontrol}) constitute a set of ordinary differential equations in terms of the states and co-states. 
The resulting TPBVP is formulated as:
\begin{align}
\boldsymbol \psi(\boldsymbol{p}_{r_0},\boldsymbol{p}_{v_0},p_{m_0}, t_f) = [\boldsymbol{r}(t_f),\boldsymbol{v}(t_f),p_m(t_f),\mathscr H(t_f)],
\label{EQ:TPBVP}
\end{align}
where $\boldsymbol \psi$ is called the shooting function, and $\boldsymbol{p}_{r_0}$, $\boldsymbol{p}_{v_0}$, and $t_f$ is the initial guess of $\boldsymbol{p}_{r}$, $\boldsymbol{p}_{v}$, and the final time, respectively.
\section{With the Final Steering Angle Constraint}\label{sec4}
Solving Eq.~(\ref{EQ:TPBVP}) will provide only fuel optimal solution with the constraint in Eq.~(\ref{EQ:steeringcondition}) unsatisfied. To this end, we modify the cost functional in Eq.~(\ref{EQ:cost}) by augmenting a regularization term defined as:
\begin{align}
\Delta (z,\theta,\beta,\epsilon) :=  \frac{1}{2}\exp(\beta z)\frac{\theta^2}{z+\epsilon},
\label{EQ:new_cost_term}
\end{align}
where $\beta$ is a small negative constant and $\epsilon$ is a very small positive constant to avoid numerical singularity. In this way, the augmented cost functional $\tilde{J}$ is
\begin{align}
\tilde{J} = \int_0^{\tilde{t}_f} (1+\Delta) u(\tilde{t}) ~\mathrm{d}\tilde{t} = \int_0^{\tilde{t}_f} \left\{ 1 + \frac{1}{2}\exp[\beta \tilde{z}(\tilde{t})]\frac{\tilde{\theta}^2(\tilde{t})}{\tilde{z}(\tilde{t})+\epsilon}\right\} \tilde{u}(\tilde{t}) ~\mathrm{d}\tilde{t}, 
\label{EQ:new_cost}
\end{align}
in which the symbol $\tilde~$ will be used in the context of the augmented cost functional hereafter.

Now, we shall show that the following property holds for the regularization term $\Delta(\tilde{z}(\tilde{t}),\tilde{\theta}(\tilde{t}),\beta,\epsilon)$: 

{\it Property:} The term $\Delta(\tilde{z}(\tilde{t}),\tilde{\theta}(\tilde{t}),\beta,\epsilon)$ remains nonnegative and sufficiently small 
for $\tilde{t} \in [0, \tilde{t}_f]$.

{\it Proof:} Notice that $\tilde{z}(\tilde{t}) > 0$ holds before touchdown and $\beta$ is a small negative constant, then the term $\exp[\beta \tilde{z}(\tilde{t})]$ becomes a very small positive number. It is clear that $\frac{\tilde{\theta}^2(\tilde{t})}{\tilde{z}(\tilde{t})+\epsilon} \geq 0$ always holds during the entire landing.
To minimize the cost functional $\tilde{J}$, the nonnegative term $\frac{\tilde{\theta}^2(\tilde{t})}{\tilde{z}(\tilde{t})+\epsilon}$ cannot approach infinity.
At touchdown $\tilde{t} = \tilde{t}_f$, we have $\tilde{z}(\tilde{t})=0$. In this case,  both $\lim\limits_{\tilde{t}\to \tilde{t}_f} \tilde{\theta}(\tilde{t}) = 0$ and  $\lim\limits_{\tilde{t}\to \tilde{t}_f} \frac{\tilde{\theta}^2(\tilde{t})}{\tilde{z}(\tilde{t})+\epsilon} = 0$ have to hold in order for the term $\Delta(\tilde{z}(\tilde{t}),\tilde{\theta}(\tilde{t}),\beta,\epsilon)$ not to blow up. This completes the proof. $\square$

Thus, the final steering angle constraint in Eq.~(\ref{EQ:steeringcondition}) will be met automatically without causing too much penalty on the original cost functional in Eq.~(\ref{EQ:cost}). In what follows, we will derive the corresponding TPBVP concerning the augmented cost functional $\tilde{J}$.

Notice that the dynamics of the lunar lander remains the same as in Eq.~(\ref{EQ:dynamics}), i.e.,
\begin{align}
\dot{\tilde{\boldsymbol x}}(\tilde t) = \boldsymbol f(\tilde{\boldsymbol x}, \tilde u, \tilde{\theta},\tilde t).
\label{EQ:new_dynamics}
\end{align}
Then, the augmented Hamiltonian, denoted by $\tilde{\mathscr{H}}$, is expressed as
\begin{align}
\tilde{\mathscr{H}} = \tilde{\boldsymbol{p}}^T_r \cdot \tilde{\boldsymbol{v}}  + \tilde{\boldsymbol{p}}^T_v \cdot (\boldsymbol g + \frac{\tilde{u}T_{m}}{\tilde{m}}\hat{\boldsymbol i}_{\tilde{\theta}}) + \tilde{p}_{m} (-\frac{\tilde{u}T_{m}}{I_{sp}g_0}) + [ 1 + \frac{1}{2}\exp(\beta \tilde{z})\frac{\tilde{\theta}^2}{\tilde{z}+\epsilon}] \tilde{u}
\label{EQ:Ham_modified}
\end{align}
The dynamics of co-states is rebuilt as 
\begin{align}
\begin{cases}
\dot{\tilde{p}}_y(\tilde t) = -\frac{\partial \tilde{\mathscr{H}}}{\partial \tilde{{y}}(\tilde t)} = 0,\\
\dot{\tilde{p}}_z(\tilde t) = -\frac{\partial \tilde{\mathscr{H}}}{\partial \tilde{{z}}(\tilde t)} = -\frac{1}{2} \tilde{u}(\tilde t) \tilde{\theta}^2(\tilde t)\exp[\beta\tilde{z}(\tilde t)]\frac{\beta[\tilde{z}(\tilde t) + \epsilon]-1}{[\tilde{z}(\tilde t)+\epsilon]^2},\\
\dot{\tilde{\boldsymbol{p}}}_v(\tilde t) = -\frac{\partial \tilde{\mathscr{H}}}{\partial \tilde{\boldsymbol{v}}(\tilde t)} = -\tilde{\boldsymbol{p}}_r(\tilde t),\\
\dot{\tilde{p}}_m(\tilde t) = -\frac{\partial \tilde{\mathscr{H}}}{\partial \tilde{m}(\tilde t)} =  \frac{\tilde u(\tilde t)T_m}{\tilde{m}^2(\tilde t)}\tilde{\boldsymbol{p}}^T_v(\tilde t) \cdot \hat{\boldsymbol i}_{\tilde{\theta}}(\tilde t).
\end{cases}
\label{EQ:dot_p_modified}
\end{align}
The optimal steering angle follows that
\begin{align}
\frac{\partial \tilde{\mathscr{H}}}{\partial \tilde{{\theta}}(\tilde t)} = 0.
\label{EQ:dH/dtheta_modi}
\end{align}
Explicitly rewriting Eq.~(\ref{EQ:dH/dtheta_modi}) leads to
\begin{align}
 \frac{\tilde u(\tilde t)T_m}{\tilde{m}(\tilde t)}[\tilde{p}_{v_y}(\tilde t)\cos\tilde{{\theta}}(\tilde t) - \tilde{p}_{v_z}(\tilde t)\sin\tilde{{\theta}}(\tilde t)] + \tilde u(\tilde t) \exp[\beta\tilde{z}(\tilde t)]\frac{\tilde{{\theta}}(\tilde t)}{\tilde{z}(\tilde t)+\epsilon} = 0.
\label{EQ:dH/dtheta_modi_real}
\end{align}
Since $\tilde u(\tilde t) \geq 0$ during landing, in order for Eq.~(\ref{EQ:dH/dtheta_modi_real}) to hold true for $\tilde t \in [0,\tilde{t}_f]$, we have
\begin{align}
 \frac{T_m}{\tilde{m}(\tilde t)}[\tilde{p}_{v_y}(\tilde t)\cos\tilde{{\theta}}^*(\tilde t) - \tilde{p}_{v_z}(\tilde t)\sin\tilde{{\theta}}^*(\tilde t)] + \exp[\beta\tilde{z}(\tilde t)] \frac{\tilde{{\theta}}^*(\tilde t)}{\tilde{z}(\tilde t)+\epsilon} = 0.
\label{EQ:dH/dtheta_modi_real_real}
\end{align}
Unlike the analytical solution in terms of $\boldsymbol{p}_v(t)$ in Eq.~(\ref{EQ:dH/dtheta}), the optimal steering angle considered here has to be determined by solving Eq.~(\ref{EQ:dH/dtheta_modi_real_real}). However, Eq.~(\ref{EQ:dH/dtheta_modi_real_real}) is a transcendental equation that may have multiple zeros. Thus, the typical Newton-like iterative method or bisection method may not find the desired zero. Based on the transforming procedure proposed by Zheng {\it et al} in Ref. \cite{zheng2021time}, we now show how to find the right zero for Eq.~(\ref{EQ:dH/dtheta_modi_real_real}).

Differentiating Eq.~(\ref{EQ:dH/dtheta_modi_real_real}) w.r.t.  $\tilde \theta$ leads to
\begin{align}
 \frac{T_m}{\tilde{m}}(-\tilde{p}_{v_y}\sin\tilde{{\theta}}^* - \tilde{p}_{v_z}\cos\tilde{{\theta}}^*) + \exp(\beta\tilde{z}) \frac{1}{\tilde{z}+\epsilon} = 0.
\label{EQ:dH/dtheta_modi_real_zeros}
\end{align}
By substituting the half-angle formulas
\begin{align*}
\sin\tilde{{\theta}} = \frac{2\tan\frac{\tilde{{\theta}}}{2}}{1+\tan^2\frac{\tilde{{\theta}}}{2}} ~~ {\rm and} ~~
\cos\tilde{{\theta}} = \frac{1-\tan^2\frac{\tilde{{\theta}}}{2}}{1+\tan^2\frac{\tilde{{\theta}}}{2}},
\end{align*}
into Eq.~(\ref{EQ:dH/dtheta_modi_real_zeros}), one has that $\tan\frac{\tilde{{\theta}}}{2}$ is a zero of the following quadratic polynomial in terms of $x$:
\begin{align}
[-\tilde{p}_{v_z} - \frac{\tilde{m}\exp(\beta\tilde{z})}{T_m(\tilde{z}+\epsilon)}]x^2+2\tilde{p}_{v_y}x + \tilde{p}_{v_z} - \frac{\tilde{m}\exp(\beta\tilde{z})}{T_m(\tilde{z}+\epsilon)} = 0.
\label{EQ:dH/dtheta_modi_real_quad}
\end{align}
It is evident that the two roots of Eq.~(\ref{EQ:dH/dtheta_modi_real_quad}) can be obtained either by radicals or standard polynomial solvers. A bisection method can be then used to find all the zeros of Eq.~(\ref{EQ:dH/dtheta_modi_real_real}) by comparing the values for Eq.~(\ref{EQ:dH/dtheta_modi_real_real}) when $\tilde{{\theta}}$ takes the values of $-\pi$, $\pi$, and the two roots of Eq.~(\ref{EQ:dH/dtheta_modi_real_quad}). After all the zeros of Eq.~(\ref{EQ:dH/dtheta_modi_real_real}) are found, the optimal value for $\tilde{{\theta}}$ is whichever that minimizes the  Hamiltonian in Eq.~(\ref{EQ:Ham_modified}).

Since the thrust ratio is still linear w.r.t. both the state dynamics in Eq.~(\ref{EQ:dynamics}) and the augmented cost functional in Eq.~(\ref{EQ:new_cost}), it is reasonable to assume that the optimal thrust ratio is still in a bang-bang form, i.e.,
\begin{align}
\tilde{u}^*(\tilde t) = 
\begin{cases}
1, \tilde{S}(\tilde t) < 0, \\
0, \tilde{S}(\tilde t) > 0,
\end{cases}
\label{EQ:optimalration_modi}
\end{align}
where $\tilde{S}(\tilde t)$ is defined as
\begin{align}
\tilde{S}(\tilde t) =  \tilde{\boldsymbol{p}}^T_v \cdot \hat{\boldsymbol i}_{\tilde{\theta}} \frac{T_m}{\tilde{m}}  - \frac{T_m}{I_{sp}g_0} \tilde{p}_{m} + 1 + \frac{1}{2}\exp(\beta \tilde{z})\frac{\tilde{\theta}^2}{\tilde{z}+\epsilon}.
\label{EQ:SF_modified}
\end{align}
Likewise, the discontinuous bang-bang control in Eq.~(\ref{EQ:optimalration_modi}) is approximated by
\begin{align}
\tilde{u}^*(\tilde t) \approx \tilde{u}^*(\tilde t,\delta) = \frac{1}{2}(1 - \frac{\tilde{S}(\tilde t)}{\sqrt{\delta + |\tilde{S}(\tilde t)|^2}}),
\label{EQ:smoothcontrol2}
\end{align}
in which $\delta$ is the same smoothing constant as in Eq.~(\ref{EQ:smoothcontrol}).
The transversality condition concerning $\tilde p_{m}$ stays unchanged, i.e., 
\begin{align}
\tilde{p}_{m}(\tilde{t}_f) = 0.
\label{EQ:Transversality1_modif}
\end{align}
Meanwhile, the stationary condition is
\begin{align}
\tilde{\mathscr{H}}(\tilde{t}_f) = 0.
\label{EQ:stat1_modif}
\end{align}

The initial condition in Eq.~(\ref{EQ:initialcondition}) and the final condition in Eq.~(\ref{EQ:finalcondition}) still apply. As a result, the corresponding TPBVP is rebuilt as:
\begin{align}
\boldsymbol \psi(\tilde{\boldsymbol{p}}_{r_0},\tilde{\boldsymbol{p}}_{v_0},\tilde{p}_{m_0}, \tilde{t}_f) = [\tilde{\boldsymbol{r}}(\tilde{t}_f),\tilde{\boldsymbol{v}}(\tilde{t}_f),\tilde{p}_m(\tilde{t}_f),\tilde{\mathscr{H}}(\tilde{t}_f)],
\label{EQ:TPBVP_modified}
\end{align}
where $\tilde{\boldsymbol{p}}_{r_0}$, $\tilde{\boldsymbol{p}}_{v_0}$, $\tilde{p}_{m_0}$, and $\tilde{t}_f$ is the initial guess of $\tilde{\boldsymbol{p}}_{r}$, $\tilde{\boldsymbol{p}}_{v}$, $\tilde{p}_{m}$, and the final time, respectively.

\section{Numerical Simulations}\label{sec5}
Before proceeding, we shall present some constants for simulations. The propulsion system of the lunar lander is specified by $I_{sp} = 311$ s and $T_m = 44,000$ N. The gravitational acceleration of the Moon is $g = 1.6229$ $\rm m/s^2$, and $g_e$ is equal to 9.81 $\rm m/s^2$. The absolute and relative tolerances for propagating the dynamics of states and co-states are set as $1.0 \times 10^{-10}$. The termination tolerance for {\it fsolve} is set as $1.0 \times 10^{-10}$. The smoothing constant $\delta$ is $1.0 \times 10^{-10}$. The constants in Eq.~(\ref{EQ:new_cost}) are $\beta = 1.0 \times 10^{-2}$ and $\epsilon = 1.0 \times 10^{-8}$.

\begin{figure}[!htp]
\centering
\begin{subfigure}[t]{7cm}
\centering
\includegraphics[width = 8cm]{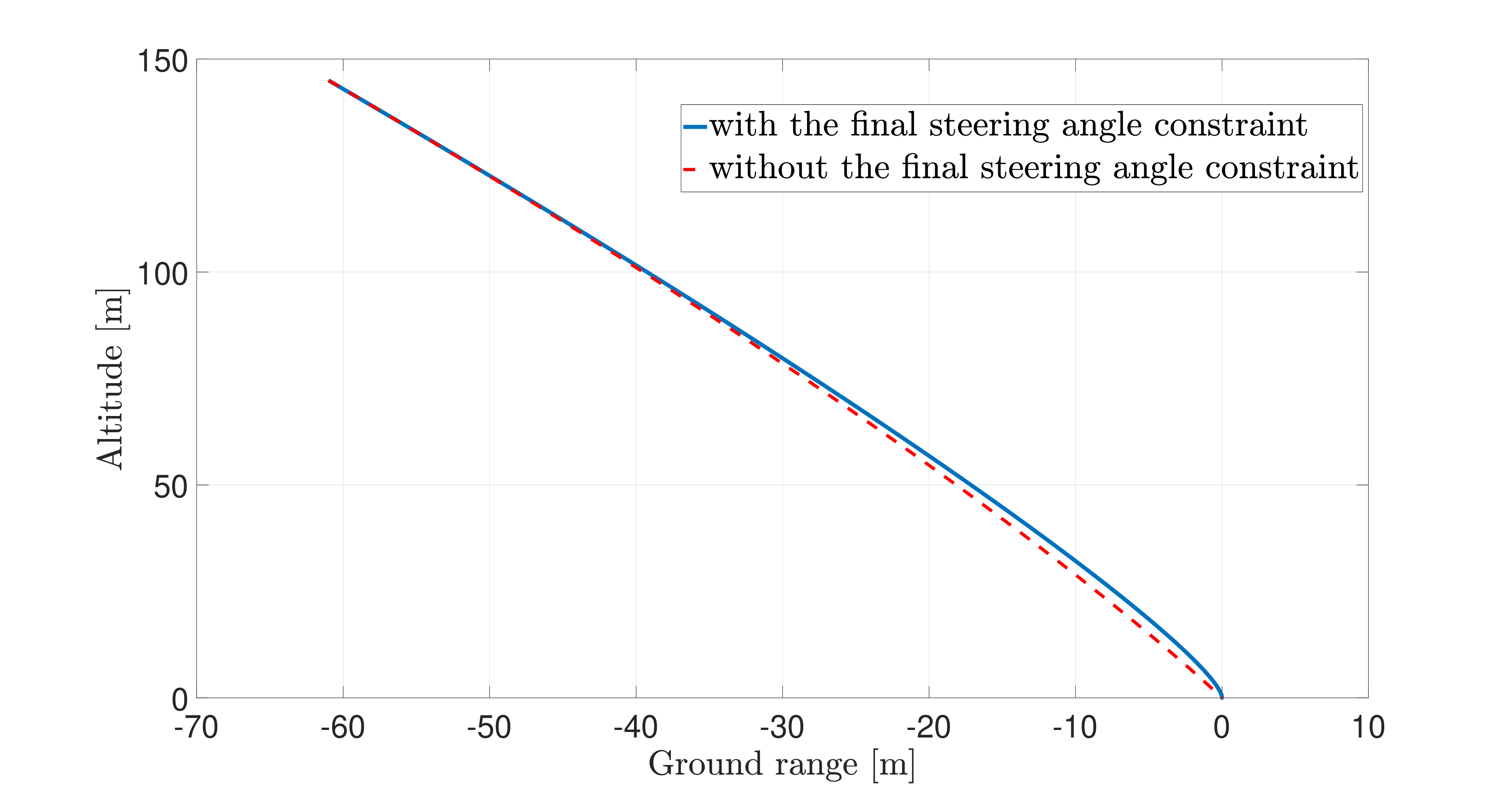}
\caption{Flight trajectories.}
\label{Fig:cooperative_control_1}
\end{subfigure}
~~~~~
\begin{subfigure}[t]{7cm}
\centering
\includegraphics[width = 8cm]{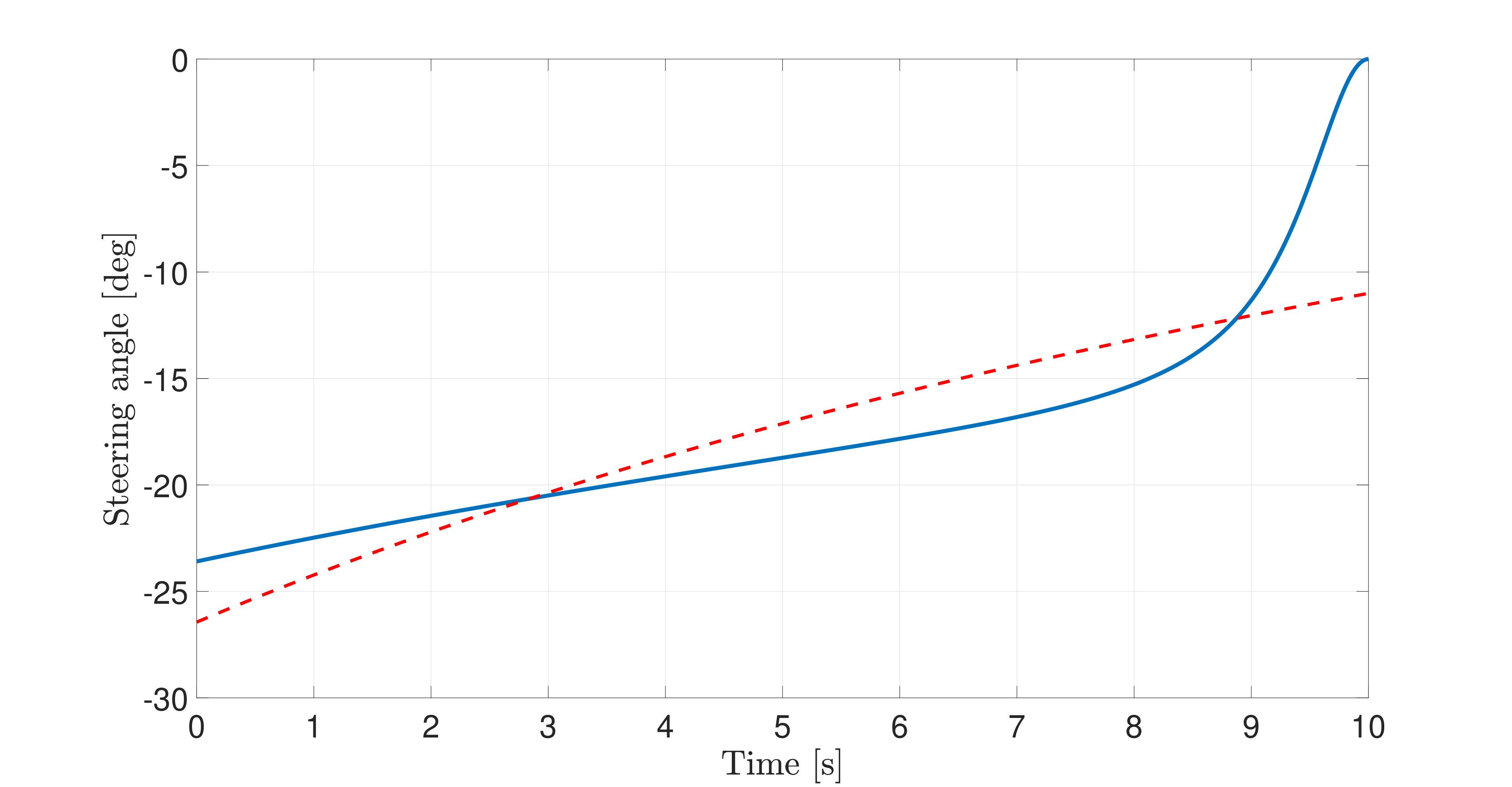}
\caption{Profiles of the steering angle.}
\label{Fig:cooperative_control_2}
\end{subfigure}\\
\begin{subfigure}[t]{7cm}
\centering
\includegraphics[width = 8cm]{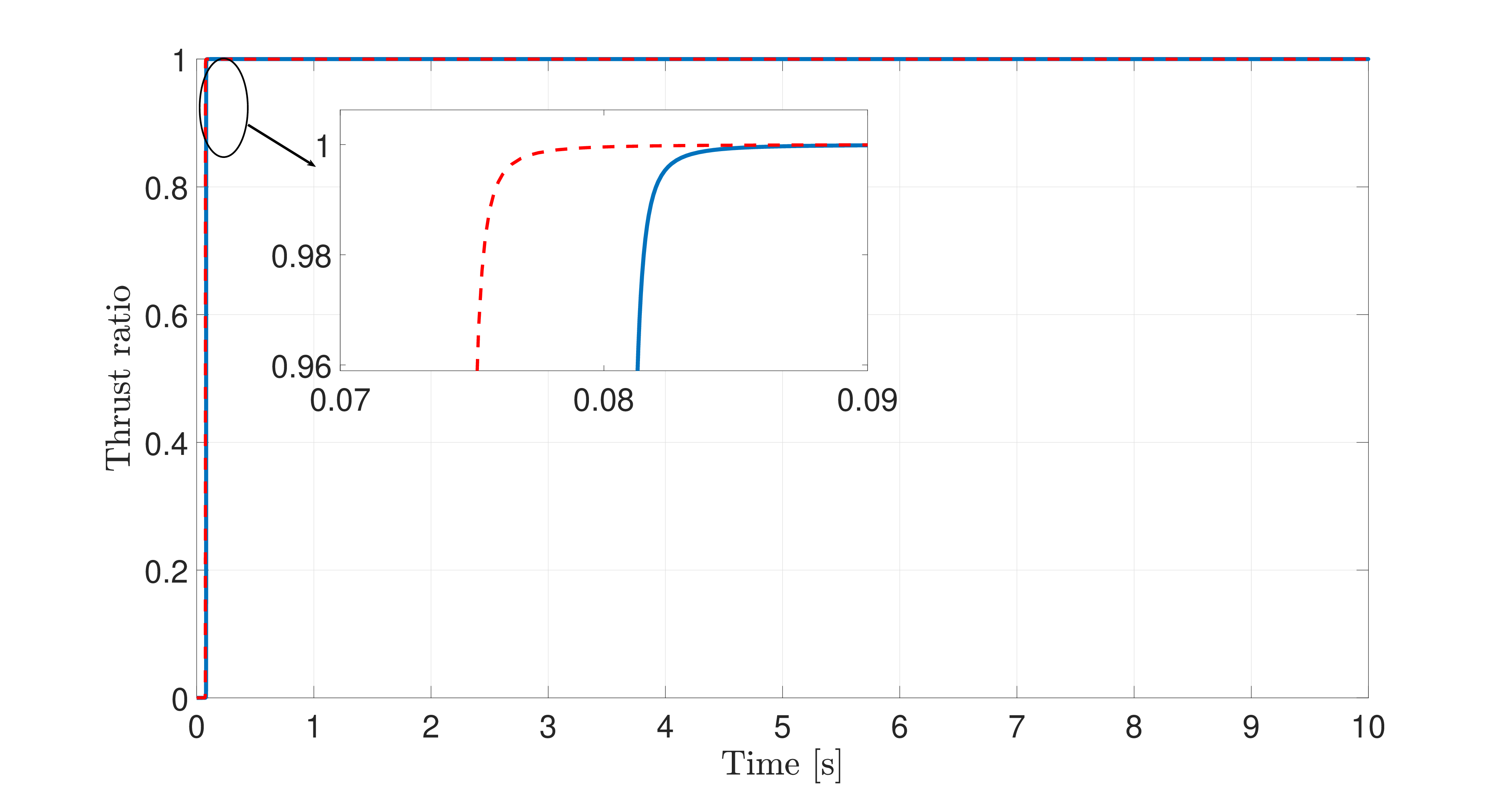}
\caption{Profiles of the thrust ratio.}
\label{Fig:cooperative_control_3}
\end{subfigure}
~~~~~
\begin{subfigure}[t]{7cm}
\centering
\includegraphics[width = 8cm]{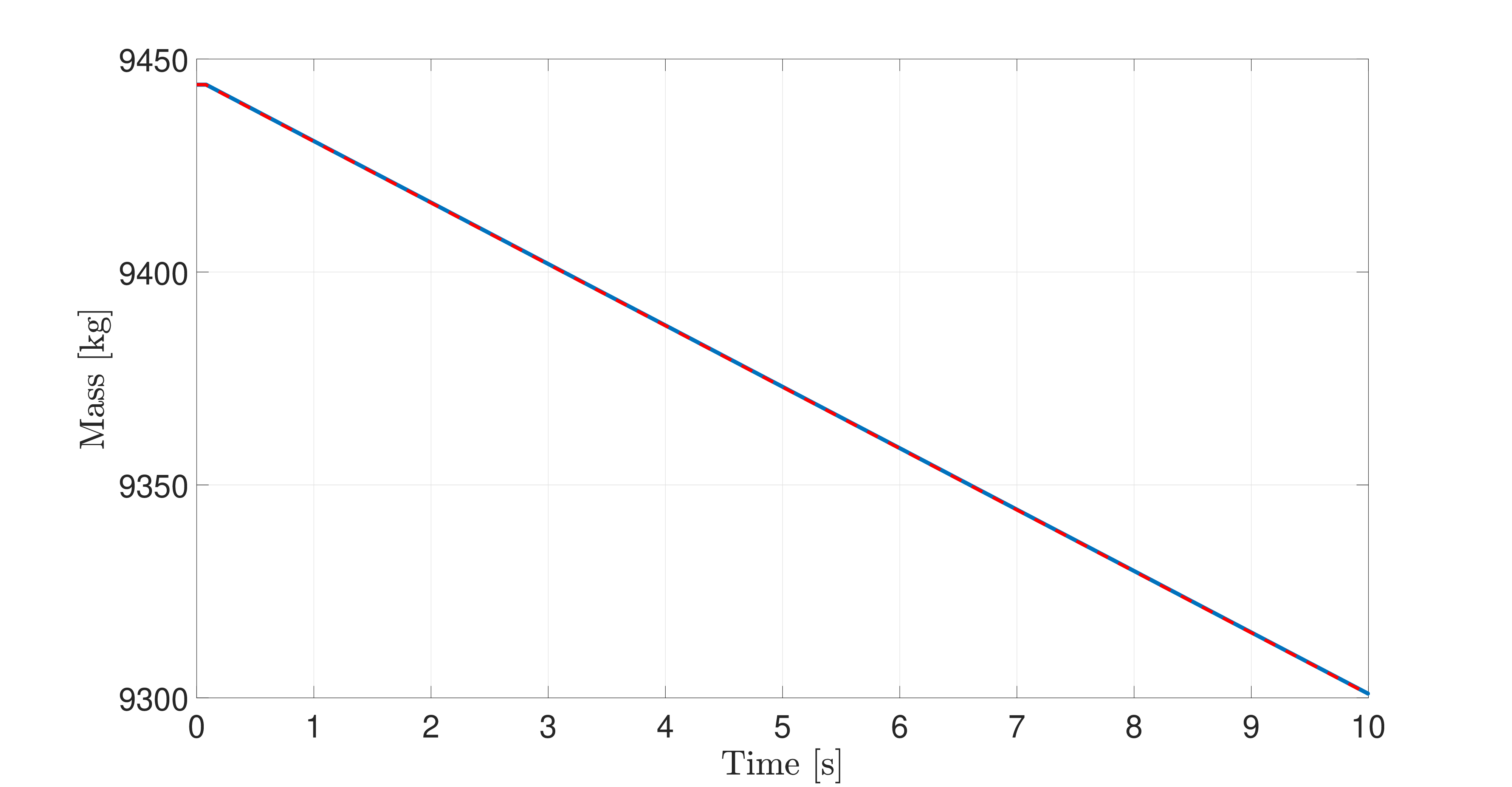}
\caption{Profiles of the mass.}
\label{Fig:cooperative_control_4}
\end{subfigure}
\caption{Trajectories, profiles of steering angle, thrust ratio, and mass with and without the final steering angle constraint.}
\label{Fig:cooperative_profile}
\end{figure}

The initial condition of the lunar lander is set as $\boldsymbol{x}_0 = [-61~{\rm m}, 145~{\rm m}, 14~{\rm m/s}, -28~{\rm m/s}, 9444~{\rm kg}]$. Solving the shooting function in Eqs.~(\ref{EQ:TPBVP}) and (\ref{EQ:TPBVP_modified}) enables one to obtain the optimal solution without and with the final steering angle constraint. Fig.~\ref{Fig:cooperative_profile} displays the flight trajectories, profiles of steering angle and thrust ratio, as well as the mass during landing. As seen from Figs.~\ref{Fig:cooperative_control_1} and \ref{Fig:cooperative_control_2}, the flight trajectories are different, in which the final steering angle obtained from Eq.~(\ref{EQ:TPBVP_modified}) reaches zero as expected; in contrast, the steering angle obtained from Eq.~(\ref{EQ:TPBVP}) is progressing linearly and arrives at $-11.02$ deg at touchdown. As for the thrust ratio, both solutions show an ``off-on'' structure. However, the switching time is slightly different; the engine is switched from off to on at $\tilde{t} = 0.0811$ s for the solution with the final steering angle constraint. On the other hand, for the solution without the final steering angle constraint, the engine is switched on at $t = 0.0748$ s, as shown by the  scaled plot in Fig.~\ref{Fig:cooperative_control_3}. 
As a result, the final time for the solution with and without the final steering angle constraint is $\tilde{t}_f = 9.9994$ s and $t_f = 9.9779$ s, respectively. Moreover, the final mass for these two solutions is $\tilde m(\tilde{t}_f) = 9300.96$ kg and $m(t_f) = 9301.18$ kg, respectively. This indicates that adding the regularization term in Eq.~(\ref{EQ:new_cost_term}) only results in an extra fuel consumption of $0.22$ kg, which is negligible in the context of an initial mass of $9444$ kg.

\begin{figure}[!htp]
\begin{center}
\includegraphics[scale=0.2]{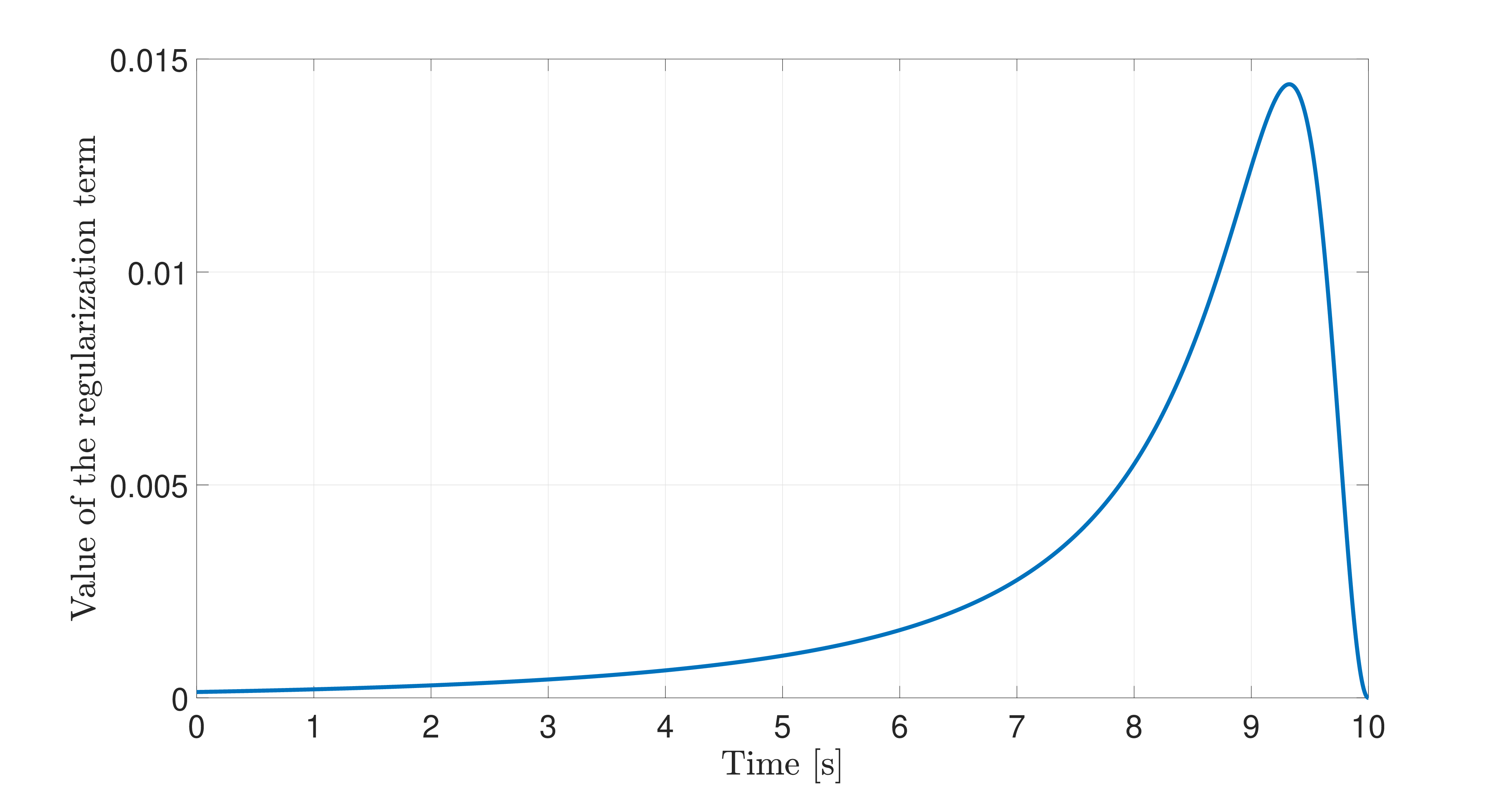}
\caption{Profile of the regularization term during landing.}\label{Fig:estimate_regularization}
\end{center}
\end{figure}

\begin{figure}[!htp]
\centering
\begin{subfigure}[t]{0.45\textwidth}
\centering
\includegraphics[scale=0.126]{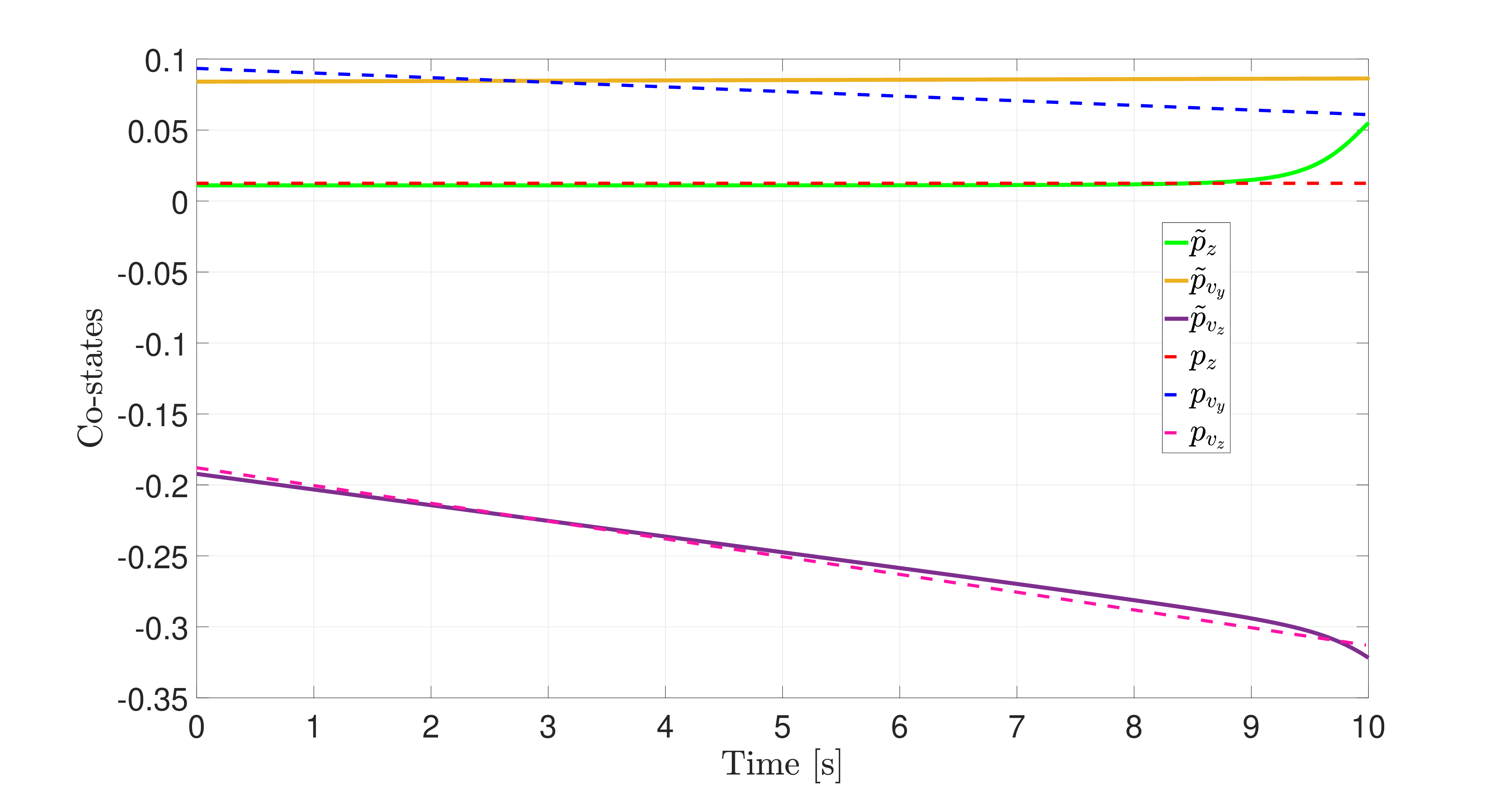}
\caption{Comparisons of some co-states.}
\label{Fig:direction_shooting}
\end{subfigure}
~~~~~
\begin{subfigure}[t]{0.45\textwidth}
\centering
\includegraphics[scale=0.126]{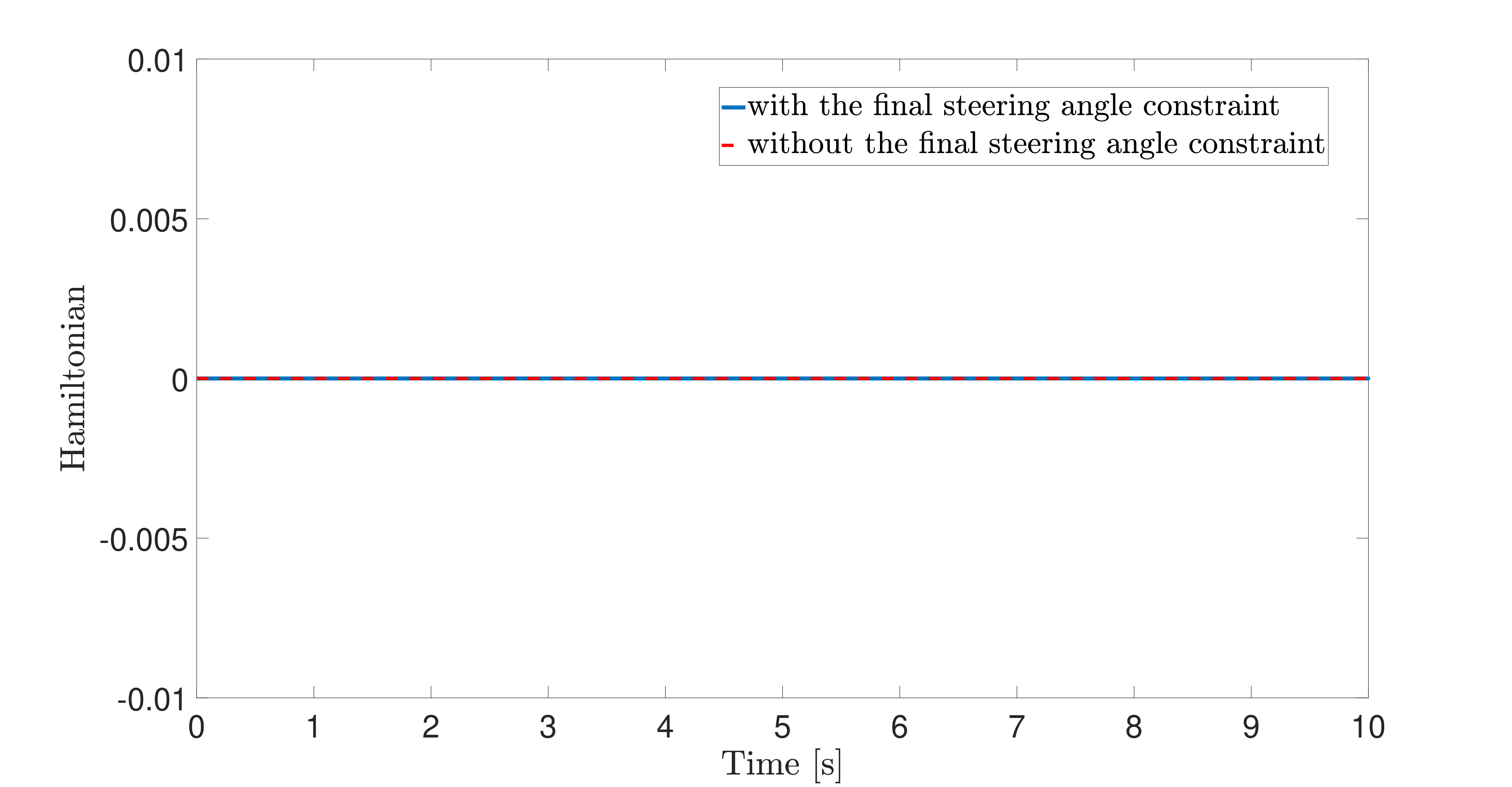}
\caption{Comparisons of the Hamiltonian.}
\label{Fig:magnitude_shooting}
\end{subfigure}\\
\caption{Comparisons of some co-states and the Hamiltonian with and without the final steering angle constraint.}
\label{Fig:ShootingSolution}
\end{figure}

In addition, Fig.~\ref{Fig:estimate_regularization} shows the profile of the regularization term defined in Eq.~(\ref{EQ:new_cost_term}), from which we can see that the regularization term $ \Delta $ remains nonnegative and small during the entire landing. Specifically, before touchdown, $\Delta$ is always a small positive number; while at touchdown, $\Delta (\tilde{t}_f)$ becomes zero because the final steering angle $\tilde{\theta}(\tilde{t}_f)$ becomes zero. 
Fig.~\ref{Fig:ShootingSolution} compares some co-states and the Hamiltonian for these two solutions. It is clear from Fig.~\ref{Fig:direction_shooting} that $p_z$ is constant, while $p_{v_y}$ and $p_{v_z}$ change linearly w.r.t time, as predicted by Eq.~(\ref{EQ:dot_p}). On the other hand, for the solution with the final steering angle constraint, both $\tilde{p}_{v_y}$ and $\tilde{p}_{v_z}$ are progressing nonlinearly w.r.t time. Moreover, the optimality of the obtained solutions can be further verified by the Hamiltonian in Fig.~\ref{Fig:magnitude_shooting}, where both solutions remain zero during the entire landing.

\begin{figure}[!htp]
\centering
\begin{subfigure}[t]{7cm}
\centering
\includegraphics[width = 8cm]{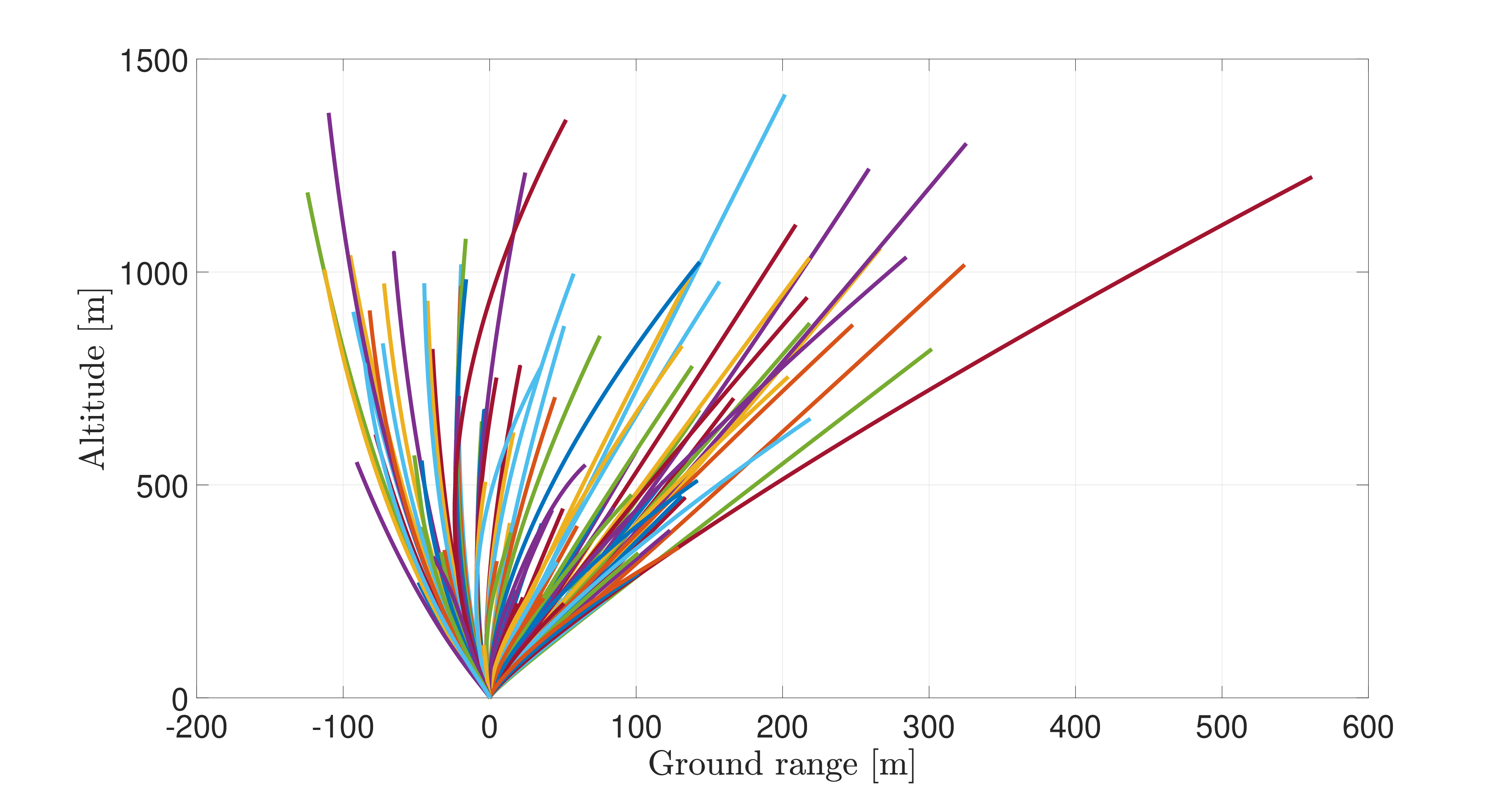}
\caption{Flight trajectories.}
\label{Fig:cooperative_appli_tra}
\end{subfigure}
~~~~~
\begin{subfigure}[t]{7cm}
\centering
\includegraphics[width = 8cm]{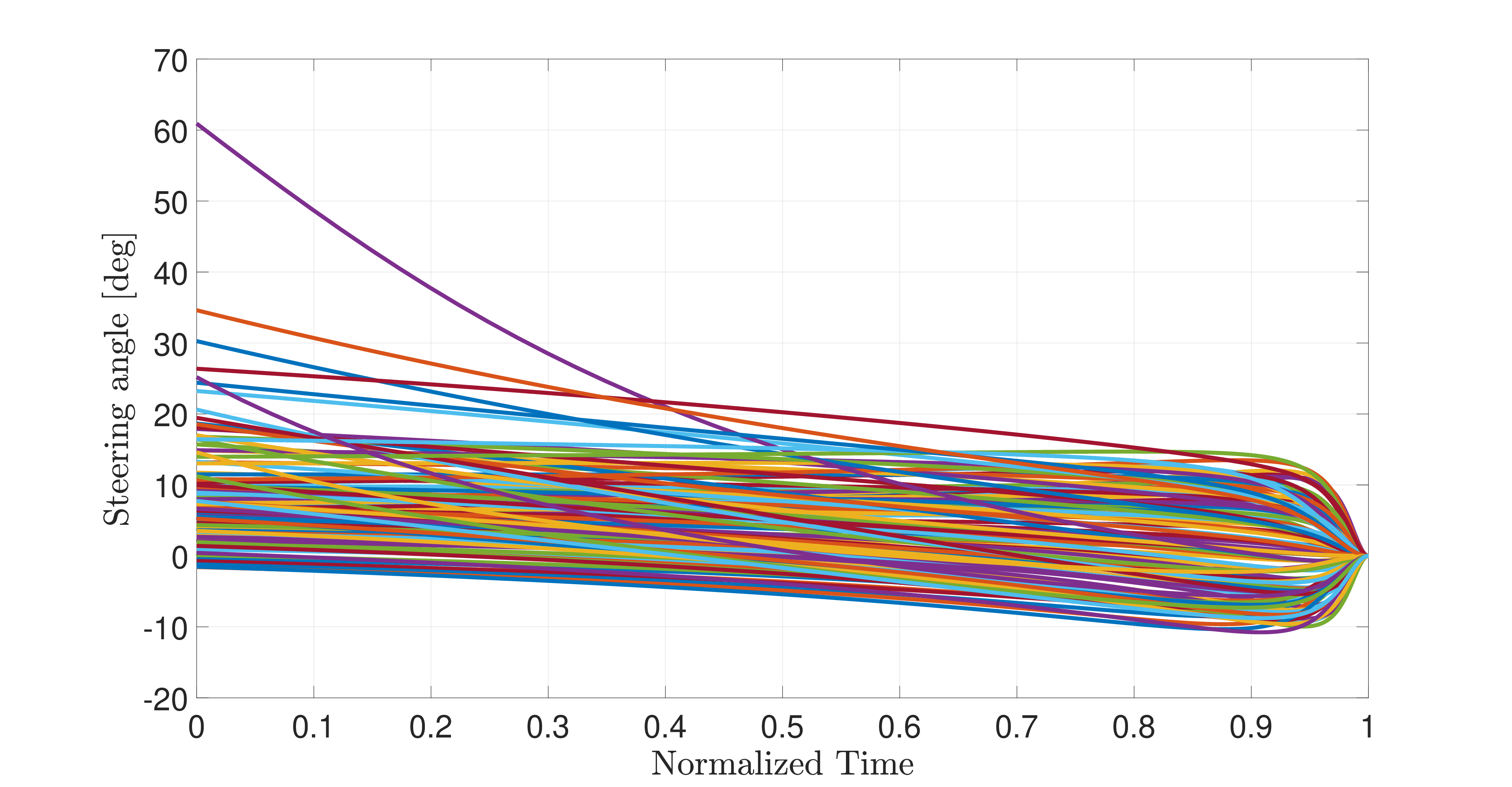}
\caption{Profiles of the steering angle.}
\label{Fig:cooperative_appli_angle}
\end{subfigure}\\
\begin{subfigure}[t]{7cm}
\centering
\includegraphics[width = 8cm]{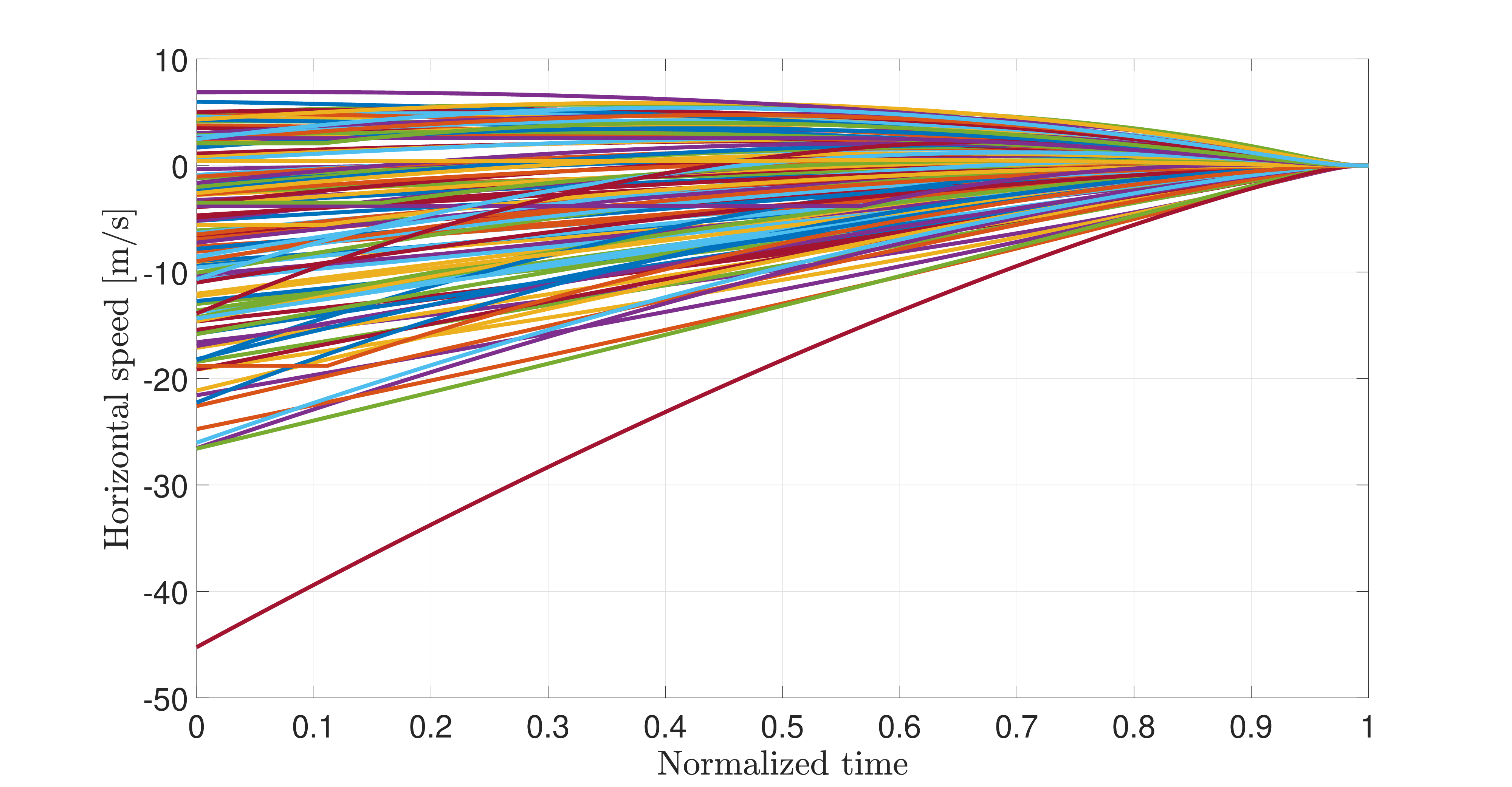}
\caption{Profiles of the horizontal speed.}
\label{Fig:cooperative_H_speed}
\end{subfigure}
~~~~~
\begin{subfigure}[t]{7cm}
\centering
\includegraphics[width = 8cm]{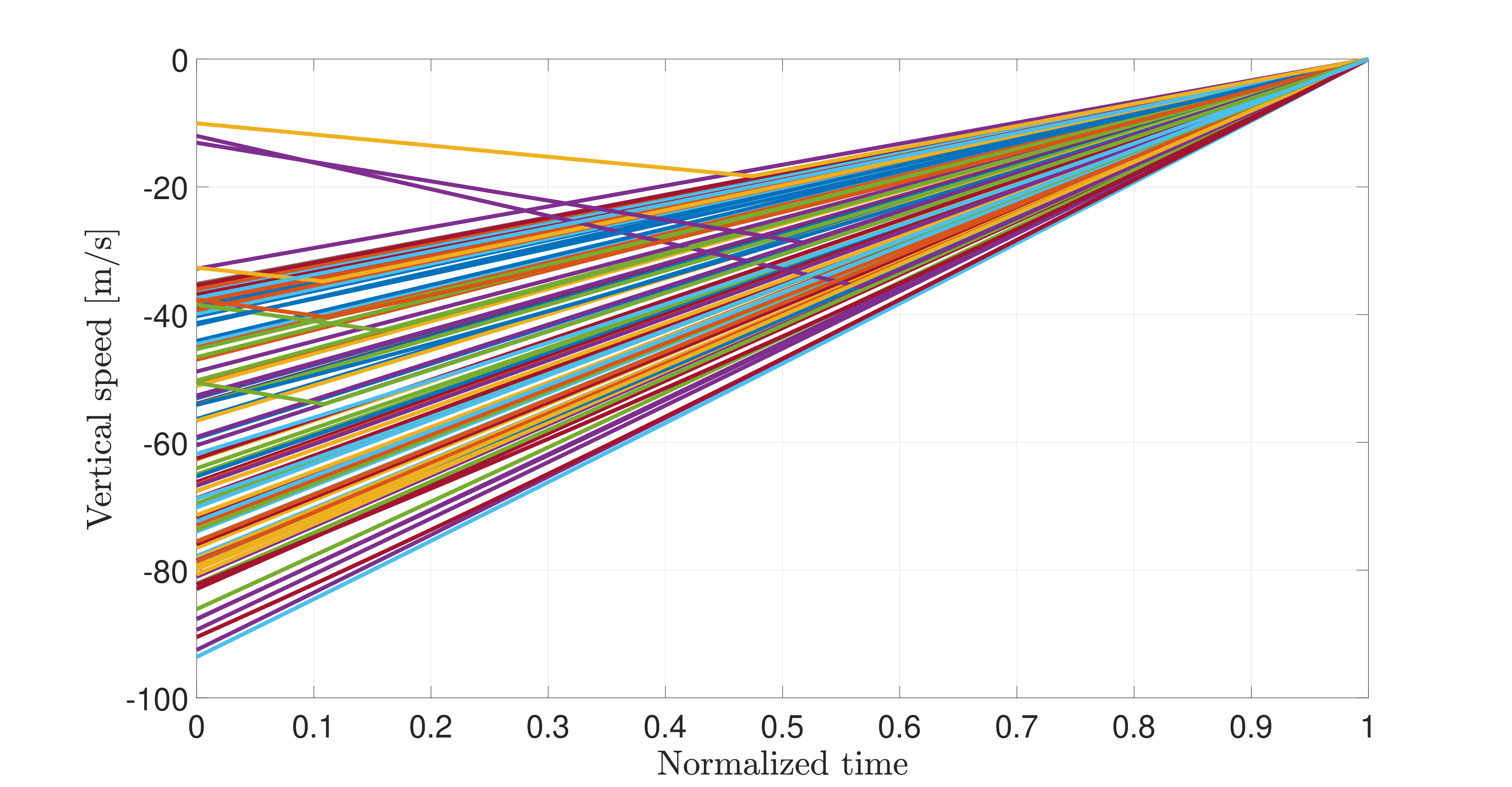}
\caption{Profiles of the vertical speed.}
\label{Fig:cooperative_V_speed}
\end{subfigure}
\caption{Trajectories, profiles of steering angle, horizontal speed, and vertical speed with the final steering angle constraint.}
\label{Fig:applicable}
\end{figure}

To further demonstrate the applicability of the proposed method, a total of 100 initial conditions are randomly generated in a domain $\mathcal{A}$, defined by 
\begin{align*}
y_0 \in [-125, 600]~({\rm m}), z_0 \in [50, 1500]~({\rm m/s}), 
v_{y_0} \in [-50, 10]~({\rm m}), v_{z_0} \in [-100, 10]~({\rm m}), m_0 \in [9050, 9450]~({\rm kg}).
\label{EQ:initial_condition_random}
\end{align*}
After solving Eq.~(\ref{EQ:TPBVP_modified}) regarding these 100 initial conditions, Fig.~\ref{Fig:applicable} shows the flight trajectories, profiles of steering angle, horizontal speed, and vertical speed. It can be observed that vertical landing is achieved for all cases, which can be further verified by the fact that the final steering angle reaches zero, as shown in
Fig.~\ref{Fig:cooperative_appli_angle}.
From Figs.~\ref{Fig:cooperative_H_speed} and \ref{Fig:cooperative_V_speed}, we can see that the horizontal speed and vertical speed become zero at touchdown. Furthermore, as for the vertical speed in most cases, it is monotonically decreasing to zero because the engine thrust is kept at the maximum during the entire landing. For these cases where the vertical speed is initially increasing solely due to gravity, the engine is switched off during this time interval. Then, the engine is kept on until touchdown.
\section{Conclusions}\label{sec6}
In this work, we proposed a method to generate the fuel-optimal vertical landing trajectory for a lunar lander. To ensure the vertical landing, we modified the cost functional by augmenting a nonnegative regularization term. In this way, in order for the Hamiltonian to be minimized, the final steering angle would  be satisfied automatically. By employing a transforming procedure, 
the necessary conditions for optimality indicated that the modified optimal steering angle could be found by a simple bisection method. As a result, the fuel-optimal vertical landing trajectory could be generated offline by the indirect shooting method. Future research directions include generalization of the proposed method to vertical landing in three dimension and its extension to onboard application.

\section*{Acknowledgement}

This research was supported by the National Natural Science Foundation of China under grant Nos. 61903331 and 62088101.

\bibliography{main}

\end{document}